# THE BRANCHING PROCESS WITH LOGISTIC GROWTH

By Amaury Lambert

*Université Pierre et Marie Curie and École Normale Supérieure, Paris*

In order to model random density-dependence in population dynamics, we construct the random analogue of the well-known logistic process in the branching process' framework. This density-dependence corresponds to intraspecific competition pressure, which is ubiquitous in ecology, and translates mathematically into a quadratic death rate. The logistic branching process, or LB-process, can thus be seen as (the mass of) a fragmentation process (corresponding to the branching mechanism) combined with constant coagulation rate (the death rate is proportional to the number of possible coalescing pairs). In the continuous state-space setting, the LB-process is a time-changed (in Lamperti's fashion) Ornstein–Uhlenbeck type process. We obtain similar results for both constructions: when natural deaths do not occur, the LB-process converges to a specified distribution; otherwise, it goes extinct a.s. In the latter case, we provide the expectation and the Laplace transform of the absorption time, as a functional of the solution of a Riccati differential equation. We also show that the quadratic regulatory term allows the LB-process to start at infinity, despite the fact that births occur infinitely often as the initial state goes to $\infty$. This result can be viewed as an extension of the pure-death process starting from infinity associated to Kingman's coalescent, when some independent fragmentation is added.

## 1. Introduction.

1.1. *Presentation of the LB-process.* The goal of this paper is to define and study the random analogue of a demographic deterministic model which is ubiquitous in ecology and widely known as the logistic growth model. It is defined by the following ordinary differential equation:

$$dZ_t = bZ_t\, dt - cZ_t^2\, dt, \qquad t > 0, \tag{1}$$

where $b$ and $c$ are real numbers ($c > 0$). The quadratic regulatory term has a deep ecological meaning as it describes negative interactions between each pair of individuals in the population. The stochastic version of the logistic process we construct is a density-dependent continuous time branching process, in both continuous and discrete state-space. We call it the "branching process with logistic growth" or "logistic branching process," abbreviated as LB-process.

In the discrete state-space setting, individuals give birth as in the pure branching case (independently, at constant rate $\rho$, to i.i.d. numbers of offspring), but can die either naturally (at constant rate $d$) or by competition pressure [when the total population has size $n$: at rate $c(n-1)$, proportionally to the number of extant conspecifics]. Therefore, when the total population has size $n$, the first birth arrives at rate $\rho n$ and the first death at rate $dn + cn(n-1)$. This size-dependence complies with the deterministic logistic growth model (1), in its formal definition (quadratic death rate), as well as in its interpretation (negative interactions between all possible pairs).

In the continuous state-space setting, the LB-process is a Markov process with (nonnegative) real values. Its paths are a.s. cdlg, that is, they are right-continuous with left-hand limits. The general definition is inspired by Lamperti's transform [27] linking continuous-state branching processes and Lévy processes (i.e., processes with independent and stationary increments), and is done by time-changing Ornstein–Uhlenbeck type processes ([34], page 104). In the special case when the LB-process has a.s. continuous paths, it is the unique strong solution (when the initial state is fixed) of a stochastic differential equation (SDE) of the following type:

$$dZ_t = bZ_t\, dt - cZ_t^2\, dt + \sqrt{\gamma Z_t}\, dB_t, \qquad t > 0, \tag{2}$$

where $B$ is the standard Brownian motion and $\gamma$ the so-called Gaussian coefficient. In this case, we retain the name coined by Etheridge in [13] where she studies it in a spatial setting: "Feller diffusion with logistic growth."

1.2. *Modeling density-dependence.* In population biology, the most simple process modeling the dynamics of a population is the Malthusian process. If $Z_t \in [0, \infty)$ denotes the total number of individuals at time $t$, then the Malthusian process satisfies $dZ_t = bZ_t\, dt$, where $b$ is the mean birth–death balance per individual and per time unit. The solutions are straightforward exponential functions and when $b > 0$, they rapidly go to $\infty$, proving useless for long-term models. Moreover, this model does not allow populations with positive growth to become extinct.

This elementary model has a random counterpart, called the branching process, where populations may have positive (expected) growth and become extinct. In their discrete time and discrete state-space form, branching processes go back to Lord Francis Galton and Irénée-Jules Bienaymé. The so-called Bienaymé–Galton–Watson (BGW) process is a Markov chain, where time steps are the nonoverlapping generations, with individuals behaving independently from one another and each giving birth to a (random) number of offspring (belonging to the next generation). These (random) offspring all have the same probability distribution. Here, the mean growth is geometric, but the process evolves randomly through time, eventually dying out or tending to $\infty$, with probability 1.

Despite its advantage of allowing extinctions while the mean growth is positive, the BGW-process shares with the Malthusian process the shortcoming (from an ecological standpoint) of being able to go to $\infty$. In the deterministic case, a celebrated improvement of the Malthusian process is the logistic process (1). It is an elementary combination of geometric growth for small population sizes and a quadratic density-dependent regulatory mechanism.

The main advantage of this model is that $Z_t$ converges to a finite limit as $t \to \infty$, namely, $b/c$ (if $b > 0$) or 0 (if $b \leq 0$). On the other hand, this model does not allow the population to evolve once it has reached its stable state. A natural continuation will then be for us to replace geometric growth in the logistic equation by random branching (random growth with geometric mean). Alternatively, this can be seen as improving the branching process by, loosely speaking, adding a quadratic regulatory term to it (and thus prevent it from going to $\infty$).

It is actually a general feature of numerous models in population dynamics that the process describing the evolution of the population size through time either goes to $\infty$ or ultimately dies out [17]. In particular, this is the case for any integer-valued Markov chain with 0 as an accessible, absorbing state, including (most) BGW processes. As a consequence, the question of adding density-dependence to branching processes is not a new concern, and answers can roughly be divided into three types (a brief survey in this matter can be found in [18]). The first approach sticks to the branching scheme in the discrete time setting: at each generation, individuals have i.i.d. random numbers of offspring, but their common distribution depends on the current population size (see, e.g., [23] and the references therein; see also [14, 35]; a special Poisson case is treated in [7] where, despite the title of his paper, the author acknowledges that the model is indeed a density-dependent discrete time branching process, but has no relation to the logistic growth model (1)). The second approach relies on stochastic difference equations, giving rise to what is called the controlled BGW process (see, e.g., [16, 21, 33]). The third approach consists in generalizing the birth and death rates of the branching process in continuous time in various ways: by considering polynomial rates as functions of the population size instead of linear rates, or by setting to zero the birth rate outside a compact set of population sizes. This way of modeling density dependence is popular among biologists (see, e.g., [11, 29, 30]).

The present work pertains to the latter approach, but here we are interested in a specific phenomenon, namely, constant pairwise competition pressure, which is a ubiquitous fundamental mechanism in biology in various space and time scales, and thus deserves special attention. We will consider branching processes in their continuous time form, in both continuous and discrete state-space (note, nevertheless, that a similar construction in discrete time could be done, by sticking to the branching scheme but letting the probability of having no child depend on the population size in the same fashion).

1.3. *Statement of results.* Because all properties of the integer-valued LB-process are common to the real-valued LB-process and not conversely, we state the results in the continuous setting. Roughly speaking, these results hold under condition (L) that the logarithm of the offspring size at each birth time has finite expectation.

Consider the LB-process with continuous paths presented above (2). We will show that this process is a time-change in Lamperti's fashion [27] of the Ornstein–Uhlenbeck process, which is itself defined from a positive real number $c$, termed the competititon rate, and the Lévy process $(\sqrt{\gamma}B_t + bt, t \geq 0)$. When the paths need not be continuous, we will consider a general Lévy process with no negative jumps instead of the previous Brownian motion with drift. The LB-process is thus characterized by $c > 0$ and a real-valued function $\psi$, termed the branching mechanism, which is the Lévy exponent of a Lévy process with no negative jumps. In the case of Feller's diffusion with logistic growth (2), $\psi(\lambda) = \frac{\gamma}{2}\lambda^2 - b\lambda$, $\lambda \geq 0$.

Let us briefly sum up the results. When the underlying branching mechanism does not allow the population to decrease in the absence of competition pressure, that is, there are no "natural deaths" ($\psi$ is then the Lévy exponent of an increasing Lévy process or subordinator), the LB-process is recurrent. It is null-recurrent (continuous setting only) if it has zero drift (no "infinitesimal births"), Lévy measure of finite mass $\rho$, further satisfying $\rho < c$, and, in that case, it converges to 0 in probability. Otherwise it is positive-recurrent and its limiting distribution is explicitly given via its Laplace transform (Theorem 3.4; Theorem 2.2 in the discrete setting). When there are "natural deaths," the LB-process goes extinct a.s. We also display formulae for its resolvent measure and its extinction time (via Laplace transform as well, Theorem 3.9).

Let us mention that in the latter case (continuous setting only), the LB-process either goes to 0 but remains positive, or is absorbed at 0 in finite time, depending solely on the branching mechanism (Theorem 3.5), that is, according to a criterion that does not involve the competition rate $c$: absorption occurs with probability 1 if $\int^\infty d\lambda/\psi(\lambda) < \infty$, with probability 0 otherwise [where the last inequality has the standard meaning that $\int_x^\infty d\lambda/\psi(\lambda)$ exists and is finite for some $x$].

Most of all, we are able to start the LB-process at $+\infty$ (Corollary 3.10, Theorem 2.3). In the discrete setting this should be compared to the total mass process in the standard coalescent starting from infinity (see Section 2.3). Specifically, we denote by $w_q$ the unique positive solution on $(0,\infty)$ vanishing at $+\infty$, of the following Riccati differential equation:

$$y' - y^2 = -qr^2,$$

where $r$ is some positive function depending on the branching mechanism $\psi$ and the competition rate $c$ (Lemma 3.8, Lemma 2.1). Namely, set

$$\theta(\lambda) = \int_0^\lambda dt \, \exp\left(\int_0^t \frac{\psi(s)}{cs} ds\right), \qquad \lambda \geq 0,$$

and $\varphi$ its inverse function. Then $r = \varphi'/\sqrt{c\varphi}$. The standard LB-process (starting from infinity) has entrance law given by

$$\mathbb{E}_\infty(\exp(-\lambda Z_\tau)) = \exp\left(-\int_0^{\theta(\lambda)} w_q(s) \, ds\right), \qquad \lambda \geq 0, q > 0,$$

where $\tau$ is an independent exponential variable with parameter $q$. In particular, the absorption time $T_a$ has Laplace transform under $\mathbb{P}_\infty$,

$$\mathbb{E}_\infty(\exp(-qT_a)) = \exp\left(-\int_0^\infty w_q(s) \, ds\right), \qquad q > 0.$$

In addition, its expectation is finite and equal to

$$\mathbb{E}_\infty(T_a) = \int_0^\infty sr^2(s) \, ds.$$

A shortcoming of our logistic branching process is obviously that it leads to much less tractable formulae than in the pure branching case. Though, we claim that it is the most natural and realistic model for random self-regulatory population dynamics. The simple definition in the discrete setting can easily be handled by biologists to do simulations or to include these dynamics into more complicated models, such as spatial or genetic ones. Indeed, in the case when no natural deaths occur, the LB-process survives with probability 1 without tending to $\infty$, and might thus provide a fruitful framework for randomly evolving large populations, as well as everlasting genealogies. Open questions are the following: distributions of independent sums of LB-processes, conditioning of the LB-process to survive and quasi-stationary distributions, underlying genealogy of the LB-process, distribution of coalescence times.... In the pure branching case the first of these questions has an obvious answer, the others have been studied in particular in [2, 6, 10, 24, 25, 26, 28, 31].

The outline of the paper is straightforward. Section 2 is concerned with the discrete setting, and is written in a way that is supposed to not scare biologists, as it contains the results geared for potential applications. Section 3 gives a more formal account of the general continuous setting, as well as further interesting subtleties that do not appear in the discrete setting (null-recurrence, extinction without absorption). Section 4 gathers the proofs of various theorems and lemmas.

**2. The discrete setting.** In this section we define and study the LB-process living on the nonnegative integers. The definition differs from that of the LB-process in the continuous setting, but both processes have very similar properties. As a consequence, we will state the results in a rather straightforward fashion, and proofs will be omitted (we leave the slight adaptations to the reader, and refer her or him to the manuscript, to be found on the author's website, for details).

2.1. *A few definitions and a key result.* Let us remind the reader that in continuous time and in the pure branching case, individuals behave independently as in the BGW-process, but give birth (rate $\rho > 0$) and die (rate $d \geq 0$) at (random exponential) independent times. The number of offspring born at each birth time is $k$ with probability $\pi_k/\rho$, where $(\pi_k, k \geq 1)$ is a sequence of nonnegative integers such that

$$\sum_{k \geq 1} \pi_k = \rho.$$

When the population size is $n$, the total birth rate is thus $\rho n$ and the natural death rate $dn$. In order to model competition pressure, we then add extra deaths at rate $cn(n-1)$ ($c > 0$), considering that each of the $n$ individuals is in competition with the $n-1$ remaining others. These deaths occur as if each particle selected another fixed particle at constant rate $c$ and killed it. Observe also that these deaths due to competition have the same kernel as the pure-death counting process associated to Kingman's coalescent [22], so that the LB-process can be viewed as (the mass of) a combination of coagulation with independent fragmentation. This viewpoint will be developed briefly in Section 2.3.

Set

$$\psi(s) = d - (\rho + d)s + \sum_{i \geq 1} \pi_i s^{i+1}, \qquad s \in [0, 1],$$

which characterizes completely the underlying branching mechanism.

The LB-process $Z = (Z_t, t \geq 0)$ associated with $\psi$ and the positive real number $c$ is thus the (minimal Feller) process with infinitesimal generator $Q = (q_{ij}, i, j \geq 0)$, where

$$q_{ij} = \begin{cases} i\pi_{j-i}, & \text{if } i \geq 1 \text{ and } j > i, \\ di + ci(i-1), & \text{if } i \geq 1 \text{ and } j = i-1, \\ -i(d + \rho + c(i-1)), & \text{if } i \geq 1 \text{ and } j = i, \\ 0, & \text{otherwise.} \end{cases}$$

From now on, we assume condition (L) is satisfied, that is,

(L) $$\sum_{i \geq 1} \pi_i \log(i) < \infty.$$

Denote by $(\bar{\pi}_i)_i$ the tail of the measure $(\pi_i)_i$, that is, $\bar{\pi}_k = \sum_{i \geq k} \pi_i$, $k \geq 1$, and note that condition (L) is equivalent to $\sum_k (\bar{\pi}_k/k) < \infty$. It is easy to see that $\psi(s) = (1-s)(d - \sum_{k \geq 1} \bar{\pi}_k s^k)$, which allows, for any $s \in (0,1]$ under condition (L), the definition

$$m(s) = \int_s^1 \frac{\psi(v)}{cv(1-v)} \, dv.$$

Standard calculations yield

(3) $$\exp(m(s)) = \beta s^{-d/c} \exp\left( \sum_{k \geq 1} \frac{\bar{\pi}_k s^k}{ck} \right), \qquad s \in (0, 1],$$

where

$$\beta = \exp\left( -\sum_{k \geq 1} \frac{\bar{\pi}_k}{ck} \right).$$

Define the nonnegative decreasing function $\theta \colon (0,1] \mapsto [0, \xi)$ by

$$\theta(s) = \int_s^1 e^{m(v)} \, dv, \qquad v \in (0, 1],$$

where $\xi \doteq \int_0^1 e^{m(v)} \, dv \in (0, +\infty]$. As a consequence of (3), we can assert that

$$\xi = \infty \iff d \geq c.$$

The mapping $\theta$ is a bijection, whose inverse on $[0, \xi)$ will be denoted by $\varphi$. In particular,

$$\varphi'(s) = -\exp(-m \circ \varphi(s)), \qquad s \in [0, \xi).$$

LEMMA 2.1. *Assume $d \neq 0$, and denote by $(\epsilon'_q)$ the following Riccati differential equation:*

$$(\epsilon'_q) \qquad y' - y^2 = -qr^2,$$

*where*

$$r(s) = \frac{|\varphi'(s)|}{\sqrt{c\varphi(s)(1-\varphi(s))}}, \qquad s \in (0,\xi).$$

*For any positive $q$, $(\epsilon'_q)$ has a unique nonnegative solution $w_q$ defined on $(0,\xi)$ and vanishing at $\xi-$. In addition, $w_q$ is positive on $(0,\xi)$, and for any $s$ sufficiently small or large, $w_q(s) < \sqrt{q}r(s)$. As a consequence, $\int_0^\xi w_q$ converges, and $w_q$ decreases initially and ultimately.*

When the coordinate process starts from $x$, its law will be denoted by $\mathbb{P}_x$.

2.2. *Results.* The behavior of the LB-process depends on whether $d$ is positive or zero. As previously, assume condition (L) holds.

When $d = 0$, recall (3) to see that the function $\exp(m)$ can be extended continuously to $[0,1]$. It is not difficult to show that one can define the probability measure $\nu$ on $\mathbb{N}$ by

$$\exp(m(s)) = \sum_{i \geq 1} \nu_i s^{i-1}, \qquad s \in [0,1].$$

THEOREM 2.2. *When $d = 0$, the LB-process is positive-recurrent in $\mathbb{N}$, and converges in distribution to the probability measure $\mu$ defined by*

$$\mu_i = \left( \sum_{j \geq 1} j^{-1} \nu_j \right)^{-1} i^{-1} \nu_i, \qquad i \geq 1.$$

*In the binary-splitting case ($\rho = \pi_1$), this limiting distribution is that of a Poisson variable of parameter $\rho/c$ conditioned on being positive*

$$\mu_i = \frac{e^{-\rho/c}}{1 - e^{-\rho/c}} \frac{(\rho/c)^i}{i!}, \qquad i \geq 1.$$

Recall the function $w_q$ defined in Lemma 2.1. Let $T_a$ denote the absorption (extinction) time.

THEOREM 2.3. *When $d > 0$, the LB-process goes extinct a.s. Moreover, the probabilities $(\mathbb{P}_x, x \geq 0)$ converge weakly, as $x \to \infty$, to the law $\mathbb{P}_\infty$ of the standard LB-process, or logistic branching process starting from infinity. Under $\mathbb{P}_\infty$, the extinction time $T_a$ has Laplace transform*

$$\mathbb{E}_\infty(\exp(-qT_a)) = \exp\left( -\int_0^\xi w_q(z)\,dz \right), \qquad q > 0,$$

*its expectation is finite and equal to*

$$\mathbb{E}_\infty(T_a) = \int_0^\xi sr^2(s)\,ds = \int_0^1 \frac{dv}{cv(1-v)} e^{-m(v)} \int_v^1 du\, e^{m(u)}.$$

2.3. *The link with fragmentation–coalescence processes.* Recall that when $\psi \equiv 0$ and $c > 0$, the LB-process is merely the pure-death process associated to Kingman's coalescent [22], where the coagulation rate per pair of objects is a constant equal to $2c$. As a consequence, the LB-process can be regarded in general as the mass of a coalescent process with independent fragmentation (corresponding to the branching mechanism): components of the partition coalesce at rate $2c$, split independently into a random number $k \geq 2$ of subcomponents at rate $\pi_{k-1}$, and may also disappear spontaneously at rate $d$ (erosion).

It is known that the pure-death process ($\rho = 0$) can be started at infinity, which stems from the fact that arrival times of coalescence events (coagulation from $i$ objects into $i-1$ objects) have summable expectations, proportional to $(i(i-1))^{-1}$. Then Theorem 2.3 asserts that it is still true even after adding fragmentation. Nevertheless, under $\mathbb{P}_\infty$, fragmentation events occur infinitely often with probability 1 as $t \searrow 0$. Indeed, the following probabilities are not summable:

$$\mathbb{P}_i(\text{a birth occurs before a death}) = \frac{\rho i}{(\rho + d)i + ci(i-1)} \sim \frac{\rho}{ci} \qquad \text{as } i \to \infty.$$

Let us mention a few recent results in this vein.

In the discrete reaction–diffusion model studied in [3], particles perform independent random walks, die spontaneously and locally coalesce or split with the same kernels as ours (in the special case of binary branching).

In [4], all exchangeable coalescence-fragmentation (EFC) processes living in the partitions of $\mathbb{N}$ are characterized. In a case of quick coalescence and slow fragmentation (Kingman's coefficient is nonzero and fragments always come in finite numbers), the mass of the EFC process has the same transition kernel as our (discrete) logistic branching process. This allows Berestycki to use Theorem 2.3 to show that these EFC processes "come down from infinity," in the sense that the number of components of the partition is finite at any positive time $t$ a.s., even though the process starts from dust (at $t = 0$, it is equal to the infinite collection of singletons). Note, however, that Theorem 2.3 was stated under the assumption that $d \neq 0$. To be convinced that the EFC process still comes down from infinity even when $d = 0$, observe that when there are $n \geq 2$ objects, the total death rate is $cn(n-1) \geq \frac{c}{2}n(n-1) + \frac{c}{2}n$, which is the total death rate associated to the LB-process with competition rate $c/2$ and death rate $c/2$.

Other works have considered size-biased fragmentation–coalescence processes (on the partitions of an interval; rates are proportional to the lengths of the subintervals), when fragmentation is binary. These partitions have infinitely many components and a ubiquitous equilibrium measure in this setup is the Poisson–Dirichlet distribution [9, 32].

2.4. *Convergence to the logistic Feller diffusion.* A convergence result is easily obtained that links the discrete and continuous settings. Namely, Feller's diffusion with logistic growth is the limit of a sequence (indexed by $n \in \mathbb{N}$) of discrete logistic branching processes. Rigorously, denote by $(Z_t^{(n)}, t \geq 0)$ the Markov process living on $n^{-1}\mathbb{N}$, started at $n^{-1}\lceil nx \rceil$, stopped at 0, and whose transition kernels are as follows. The integer $N_t^{(n)} = nZ_t^{(n)}$ is a binary-splitting ($\pi_1 = \rho$) LB-process with parameters (indexed by $n$) $\rho_n = \frac{\gamma}{2}n^2 + \lambda n$, $d_n = \frac{\gamma}{2}n^2 + \delta n$, $c_n = c$, where $c, \delta, \gamma, \lambda$, are positive constants. In other words:

1. $Z_t^{(n)}$ is incremented by $1/n$ at rate $(\frac{\gamma}{2}n + \lambda)n^2 Z_t^{(n)}$.
2. $Z_t^{(n)}$ is decremented by $1/n$ at rate $(\frac{\gamma}{2}n + \delta)n^2 Z_t^{(n)}$.
3. $Z_t^{(n)}$ is decremented by $1/n$ at rate $cn^2 Z_t^{(n)}(Z_t^{(n)} - n^{-1})$.

Let $b = \lambda - \delta$. Then standard results [20] show that as $n \to \infty$, the sequence $(Z_t^{(n)}, t \geq 0)_n$ converges weakly to the Feller diffusion with logistic growth started at $x$ and solution of (2),

$$dZ_t = bZ_t\, dt - cZ_t^2\, dt + \sqrt{\gamma Z_t}\, dB_t, \qquad t > 0.$$

## 3. The continuous setting.

3.1. *Preliminaries.* In the 1950s [19], the analogue of the BGW-process was defined in continuous time and continuous state-space (CB-process). All CB-processes are cdlg, and their jumps (if any) are a.s. positive. Those whose paths are continuous satisfy an SDE of the following type:

$$dZ_t = bZ_t\, dt + \sqrt{\gamma Z_t}\, dB_t, \qquad t > 0,$$

so that the definition of the LB-process in this case is clearly (2). When $b = 0$, the dynamics are referred to as quadratic branching. When $b = 0$ and $\gamma = 4$, $Z$ is the celebrated Feller diffusion or squared Bessel process of dimension 0.

Next we remind the reader of a celebrated result of Lamperti [27] which relates CB-processes and Lévy processes with no negative jumps. This result will allow us to give the appropriate definition of the branching process with logistic growth in a second section.

A Lévy process is a cdlg Markov process with independent and stationary increments. A Lévy process with continuous paths is the sum of a (scaled)

Brownian motion and a deterministic drift (possibly zero). Let $X$ be a real-valued Lévy process with no negative jumps. Let $T_0$ be the first hitting time of zero by $X$. Then define

$$\eta_t = \int_0^{t \wedge T_0} \frac{ds}{X_s}, \qquad t > 0,$$

and $(C_t, t \geq 0)$ its right-inverse. Lamperti's result then states that if

$$Z_t = X(C_t), \qquad t \geq 0,$$

then $Z$ is a CB-process, and

$$C_t = \int_0^t Z_s \, ds, \qquad t > 0.$$

Conversely, any CB-process $Z$ is a time-changed Lévy process: if $C$ is defined as above, and $\eta$ is the right-inverse of $C$, then $Z \circ \eta$ is a Lévy process with no negative jumps killed when it hits 0.

In the special case of quadratic branching, Lamperti's result predicts that if we time-change a Feller diffusion (i.e., a quadratic CB-process) in this fashion, we obtain (up to a multiplicative constant) a killed Brownian motion. It is indeed easy to check that if $Z_t = x + \int_0^t \sqrt{Z_s} \, dB_s$, then $Z \circ \eta$ is a local martingale with increasing process $\int_0^{\eta_t} Z_s \, ds = t$, and thus is a standard Brownian motion. Using this last argument, it is elementary to prove the following proposition concerning the LB-process.

PROPOSITION 3.1. *As previously, assume $Z$ is a Feller diffusion with logistic growth, that is, a diffusion solving* (2),

$$dZ_t = bZ_t \, dt - cZ_t^2 \, dt + \sqrt{\gamma Z_t} \, dB_t, \qquad t > 0.$$

*Define $\eta$ as the right-inverse of $C$, where $C_t = \int_0^t Z_s \, ds$. Then the process $R = Z \circ \eta$ solves the following SDE:*

$$dR_t = dX_t - cR_t \, dt, \qquad t > 0,$$

*where $X$ is a Brownian motion with drift, namely, $\gamma^{-1/2}(X_t - bt, t \geq 0)$ is the standard Brownian motion.*

*Conversely, let $R$ be the (strong) solution of the last SDE where $X$ is a Brownian motion with a possible drift, $T_0$ its first hitting time of 0 and $C$ the right-inverse of $\eta$, where*

$$\eta_t = \int_0^{t \wedge T_0} \frac{ds}{R_s}, \qquad t > 0.$$

*Then $Z = X \circ C$ is a diffusion process killed when it hits 0, and solves an SDE of type* (2).

This last proposition will enable us to extend the definition of branching process with logistic growth to any kind of branching mechanism.

3.2. *Construction of the logistic branching process.* Let $X$ denote a spectrally positive Lévy process (i.e., with no negative jumps). The branching mechanism function that interests us is the Lévy exponent $\psi$ of $X$ defined by

$$\mathbb{E}(\exp(-\lambda X_t)) = \exp(t\psi(\lambda)), \qquad t, \lambda \geq 0.$$

It is specified by the Lévy–Khinchin formula ([5], Chapter VII)

$$\psi(\lambda) = \alpha\lambda + \frac{\gamma}{2}\lambda^2 + \int_0^\infty (e^{-\lambda r} - 1 + \lambda r \mathbf{1}_{r<1})\Pi(dr), \qquad \lambda \geq 0,$$

where $\alpha$ is a real number, $\gamma$ a positive real number termed the Gaussian coefficient of $X$, and $\Pi$ is a positive $\sigma$-finite measure on $(0, \infty)$ such that $\int_0^\infty (1 \wedge r^2)\Pi(dr) < \infty$, termed the Lévy measure of $X$.

We now provide the definition of the LB-process.

DEFINITION 3.2. For any positive real number $c$ and any spectrally positive Lévy process $X$ with Lévy exponent $\psi$, we define the logistic branching process $\mathrm{LB}(\psi, c)$ associated with $\psi$ and $c$ starting from $x > 0$ as follows.

Let $R$ denote the unique strong solution, starting from $x$, of the following SDE:

(4) $$dR_t = dX_t - cR_t\,dt, \qquad t > 0,$$

$T_0$ its first hitting time of 0 and $C$ the right-inverse of $\eta$, where

$$\eta_t = \int_0^{t \wedge T_0} \frac{ds}{R_s}, \qquad t > 0.$$

Then the LB-process $Z$ is the Feller process defined by

$$Z_t = \begin{cases} R(C_t), & \text{if } 0 \leq t < \eta_\infty, \\ 0, & \text{if } \eta_\infty < \infty \text{ and } t \geq \eta_\infty. \end{cases}$$

If $A$ (resp. $Q$) denotes the infinitesimal generator of $X$ (resp. $Z$), then, for any differentiable $f$ in the domain of $A$,

$$Qf(z) = zAf(z) - cz^2 f'(z), \qquad z \geq 0.$$

Before continuing further, we prove some of the facts stated in this last definition. SDEs of type (4) are well studied, see, for example, [34], pages 104–113. They have strong solutions which are cdlg homogeneous strong Markov processes, known as Ornstein–Uhlenbeck type processes. An explicit formula for the unique strong solution of (4) with initial state $x$ is

(5) $$R_t = xe^{-ct} + \exp(-ct)\int_0^t \exp(cs)\,dX_s, \qquad t > 0.$$

By standard theory of Markov processes (see, e.g., [12]), $Z$ is then a cdlg time homogeneous strong Markov process. Furthermore, it is easily seen that the infinitesimal generator $U$ of $R$ is given, for any differentiable $f$ in the domain of $A$, by

$$Uf(z) = Af(z) - czf'(z), \qquad z \in \mathbb{R}.$$

As a consequence, for any time $t$ and initial condition $x \geq 0$, with $s = C_u$,

$$\begin{aligned}
\mathbb{E}_x\left(\int_0^t Qf(Z_u)\,du\right) &= \mathbb{E}_x\left(\int_0^t Uf(Z_u)Z_u\,du\right) \\
&= \mathbb{E}_x\left(\int_0^{C_t} Uf(R_s)\,ds\right) \\
&= \mathbb{E}_x(f(R(C_t))) - f(x) = \mathbb{E}_x(f(Z_t)) - f(x),
\end{aligned}$$

which shows, indeed, that $Q$ is the infinitesimal generator of $Z$.

3.3. *Properties of the logistic branching process.* In this section we are interested in the law of the LB-process, and particularly in its long-term behavior. This behavior depends on whether $X$ is a subordinator or not. We recall that a Lévy process with positive jumps is called a subordinator if it has increasing paths a.s. In that case, its paths have finite variation, its drift coefficient is nonnegative, and its Gaussian coefficient is zero.

In contrast to the pure branching case, LB-processes eventually go to 0 with probability 0 or 1. In reference to population biology, this probability is called the extinction probability (otherwise, absorption probability). When $X$ is a subordinator, this probability is 0, and the LB-process is recurrent. When $X$ is not a subordinator, it is 1, and we characterize the law of the extinction time. Unless otherwise specified, proofs of statements of this section are postponed to Section 4.

From now on, we assume condition (L) is satisfied, that is,

(L) $$\mathbb{E}(\log(X_1)) < \infty.$$

We recall that condition (L) is equivalent to $\int^\infty \log(r)\Pi(dr) < \infty$ [5]. Accordingly, it is easily seen that one can define

$$m(\lambda) = \int_0^\lambda \frac{\psi(s)}{cs}\,ds, \qquad \lambda \geq 0.$$

In the next statements, we consider the case when $X$ is a subordinator. We then denote by $\delta \geq 0$ its drift coefficient, so that ([5], Chapter III)

$$\psi(\lambda) = -\delta\lambda - \int_0^\infty \Pi(dr)(1 - e^{-\lambda r}), \qquad \lambda \geq 0.$$

We also introduce its Lévy tail $\bar{\Pi}$, that is, $\bar{\Pi}(y) = \int_y^\infty \Pi(dr)$, $y > 0$.

LEMMA 3.3. *Assume $X$ is a subordinator satisfying* (L). *Then $m$ can be expressed as*

$$-m(\lambda) = \frac{\delta}{c}\lambda + \int_0^\infty (1 - e^{-\lambda r})\frac{\bar{\Pi}(r)}{cr}\,dr, \qquad \lambda \geq 0,$$

*and the following equation*

$$\int_0^\infty \nu(dr)e^{-\lambda r} = \exp(m(\lambda)), \qquad \lambda \geq 0$$

*defines a unique probability measure $\nu$ on $(0, \infty)$. This probability measure is infinitely divisible and $\nu((\delta/c, \infty)) = 1$.*

REMARK. Note that $\nu$ is also self-decomposable. An r.v. $S$ is said to be self-decomposable if for any $a < 1$, there is an r.v. $S_a$ independent of $S$ such that

$$aS + S_a \stackrel{(\mathrm{d})}{=} S.$$

In particular, any self-decomposable distribution is infinitely divisible. Conversely, a real infinitely divisible distribution is self-decomposable if its Lévy measure is of type $|x|^{-1}k(x)\,dx$, where $k$ is increasing on $(-\infty, 0)$ and decreasing on $(0, \infty)$.

We introduce condition $(\partial)$, where $\rho$ is defined as

$$\rho \doteq \int_0^\infty \Pi(dr) \leq \infty.$$

We thus say that $(\partial)$ holds iff (at least) one of the following holds:

- $\delta \neq 0$,
- $\rho = \infty$,
- $c < \rho < \infty$.

THEOREM 3.4. *Assume $X$ is a subordinator satisfying* (L).

(i) *Assume $(\partial)$. Then the probability $\nu$ defined in Lemma* 3.3 *has $\int_0^\infty r^{-1} \times \nu(dr) < \infty$. The process $LB(\psi, c)$ is positive-recurrent in $(\delta/c, \infty)$ and converges in distribution to the probability measure $\mu$ whose size-biased distribution is $\nu$, that is,*

$$\mu(dr) = \left(\int_{(\delta/c, \infty)} s^{-1}\nu(ds)\right)^{-1} r^{-1}\nu(dr), \qquad r > 0.$$

*In particular, the expectation of the stationary probability is*

$$\int_{(\delta/c, \infty)} r\mu(dr) = \left(\int_{(\delta/c, \infty)} s^{-1}\nu(ds)\right)^{-1} < \infty.$$

(ii) *Assume (∂) does not hold. Then the process $LB(\psi,c)$ is null-recurrent in $(0,\infty)$ and converges to $0$ in probability.*

From now on, $X$ is assumed not to be a subordinator. In the next theorem, we claim that the extinction probability is 1. We set a criterion to establish whether the process remains positive or is absorbed, that is, reaches 0 in finite time. Note that the criterion for absorption does not depend on $c$ and is the same as for the branching process $(c=0)$ [15].

THEOREM 3.5. *Assume $X$ satisfies (L) and is not a subordinator. Then the process $LB(\psi,c)$ goes to $0$ a.s., and if $T_a$ denotes the absorption time*

$$T_a = \inf\{t \geq 0 : Z_t = 0\},$$

*then $\mathbb{P}(T_a < \infty) = 1$ or $0$ according to whether $\int^\infty d\lambda/\psi(\lambda)$ converges or diverges.*

The last two theorems will be proved in Section 4.

Next we intend to give deeper insight into the law of the LB-process. We will simultaneously be able to define the standard LB-process or LB-process starting from infinity.

Define $\theta : [0, +\infty) \mapsto [0, +\infty)$ by

$$\theta(\lambda) = \int_0^\lambda e^{m(s)} \, ds, \qquad \lambda \geq 0.$$

Recall [5] that when $X$ is not a subordinator, $\liminf_{\lambda \to \infty} \lambda^{-1} \psi(\lambda) > 0$. As a consequence, the mapping $\theta$ is a strictly increasing bijection, whose inverse will be denoted by $\varphi$. In particular,

$$\varphi'(\lambda) = \exp(-m \circ \varphi(\lambda)), \qquad \lambda \geq 0.$$

Define also for any nonnegative $x$ and positive $q$, the Laplace transform $G_{q,x}$ of the $q$-resolvent measure of the LB-process $(Z_t, t \geq 0)$ started at $x$,

$$G_{q,x}(\lambda) = \int_0^\infty dt \, e^{-qt} \mathbb{E}_x(e^{-\lambda Z_t}), \qquad \lambda \geq 0.$$

Note that from Definition 3.2, if $A$ still denotes the infinitesimal generator of the Lévy process $X$, then for any differentiable real function in the domain of $A$,

$$\begin{aligned}
(6) \quad & q \int_0^\infty dt \, e^{-qt} \mathbb{E}_x(f(Z_t)) \\
& \quad = f(x) + \int_0^\infty dt \, e^{-qt} \mathbb{E}_x(Z_t A f(Z_t) - c Z_t^2 f'(Z_t)).
\end{aligned}$$

Before stating the main theorem of this section, we display three lemmas, the last of which is the key result for the theorem. The first two help understand where the third one comes from. These lemmas are also proved in Section 4 (except Lemma 3.8 whose proof is replaced by the more general proof of Lemma 2.1).

LEMMA 3.6. *As a function of $\lambda \in (0, \infty)$, $G_{q,x}$ is twice continuously differentiable and solves the second-order linear differential equation $(E_q)$, where $y$ is the unknown function and $\lambda$ is the scalar variable*

$$(E_q) \qquad -c\lambda y'' + \psi y' + qy = e^{-x\lambda}.$$

LEMMA 3.7. *If $(E_q^{(h)})$ is the homogeneous differential equation associated to $(E_q)$, then:*

(i) *For any solution $(I, f_q)$ of $(E_q^{(h)})$, for any open subinterval $I_0$ of $I$ on which $f_q$ does not vanish, the function $g_q$ defined as*

$$g_q(\lambda) = -\frac{f_q'(\lambda)}{f_q(\lambda)}, \qquad \lambda \in I_0$$

*solves on $I_0$ the Riccati differential equation $(\epsilon_q)$*

$$(\epsilon_q) \qquad -y' + y^2 + \frac{\psi}{c\lambda}y = \frac{q}{c\lambda}.$$

(ii) *For any solution $(J, g_q)$ of $(\epsilon_q)$, the function $h_q$ defined as*

$$h_q(\lambda) = e^{-m \circ \varphi(\lambda)} g_q \circ \varphi(\lambda), \qquad \lambda \in \theta(J)$$

*solves on $\theta(J)$ the Riccati differential equation $(\epsilon_q')$*

$$(\epsilon_q') \qquad y' - y^2 = -qr^2,$$

*where*

$$r(\lambda) = \frac{\varphi'(\lambda)}{\sqrt{c\varphi(\lambda)}}, \qquad \lambda > 0.$$

LEMMA 3.8. *For any positive $q$, the Riccati differential equation $(\epsilon_q')$ has a unique nonnegative solution $w_q$ defined on $(0, \infty)$ and vanishing at $\infty$. In addition, $w_q$ is positive on $(0, \infty)$, and for any $\lambda$ sufficiently small or large, $w_q(\lambda) < \sqrt{q} r(\lambda)$. As a consequence, $\int_0 w_q$ converges, and $w_q$ decreases initially and ultimately.*

Now we are able to state the main results of this section.

THEOREM 3.9. *Recall the function $w_q$ defined in Lemma 3.8. For any nonnegative $\lambda$ and positive $q$, an expression for $G_{q,x}(\lambda)$ is given by*

$$qG_{q,x}(\lambda) = 1 - \int_0^{\theta(\lambda)} dt\, e^{-\int_t^{\theta(\lambda)} w_q(z)\,dz} \tag{7}$$
$$\times \int_t^\infty ds\, qr^2(s)(1 - e^{-x\varphi(s)})e^{-\int_t^s w_q(z)\,dz}.$$

*In particular, if $\int^\infty d\lambda/\psi(\lambda)$ converges, then the absorption time $T_a$ is a.s. finite and has Laplace transform under $\mathbb{P}_x$,*

$$\mathbb{E}_x(\exp(-qT_a))$$
$$= 1 - \int_0^\infty dt\, e^{-\int_t^\infty w_q(z)\,dz} \tag{8}$$
$$\times \int_t^\infty ds\, qr^2(s)(1 - e^{-x\varphi(s)})e^{-\int_t^s w_q(z)\,dz}, \qquad q > 0.$$

*In addition, its expectation is finite and is equal to*

$$\mathbb{E}_x(T_a) = \int_0^\infty ds\, sr^2(s)(1 - e^{-x\varphi(s)}) \tag{9}$$
$$= \int_0^\infty \frac{dt}{ct}(1 - e^{-tx})e^{-m(t)}\int_0^t ds\, e^{m(s)}.$$

COROLLARY 3.10. *The probabilities $(\mathbb{P}_x, x \geq 0)$ converge weakly, as $x \to \infty$, to the law $\mathbb{P}_\infty$ of the standard LB-process or logistic branching process starting from infinity. The standard LB-process is the Markov process that has the same transition kernels as the LB-process considered previously, and entrance law given by*

$$\mathbb{E}_\infty(\exp(-\lambda Z_\tau)) = \exp\left(-\int_0^{\theta(\lambda)} w_q(s)\,ds\right), \qquad \lambda \geq 0, q > 0, \tag{10}$$

*where $\tau$ is an independent exponential variable with parameter $q$. In particular, if $\int^\infty d\lambda/\psi(\lambda)$ converges, then under $\mathbb{P}_\infty$ the absorption time $T_a$ is a.s. finite and has Laplace transform*

$$\mathbb{E}_\infty(\exp(-qT_a)) = \exp\left(-\int_0^\infty w_q(s)\,ds\right), \qquad q > 0. \tag{11}$$

*In addition, its expectation is finite and equal to*

$$\mathbb{E}_\infty(T_a) = \int_0^\infty sr^2(s)\,ds = \int_0^\infty \frac{dt}{ct}e^{-m(t)}\int_0^t ds\, e^{m(s)}. \tag{12}$$

PROOF OF THEOREM 3.9. As in Lemma 3.7, define $g_q$ as

$$g_q(s) = e^{m(s)} w_q \circ \theta(s), \qquad s \geq 0.$$

Then the function $g_q$ solves the Riccati differential equation ($\epsilon_q$) and by Lemma 3.8 is integrable at 0, for

$$\int_0^\lambda g_q(s)\,ds = \int_0^{\theta(\lambda)} w_q(t)\,dt, \qquad \lambda \geq 0,$$

by the change $t = \theta(s)$. Then define the functions $F_{q,x}$ and $K_{q,x}$ as

$$F_{q,x}(\lambda) = G_{q,x}(\lambda) \exp\left(\int_0^\lambda g_q(s)\,ds\right), \qquad \lambda \geq 0,$$

$$K_{q,x}(\lambda) = F'_{q,x}(\lambda) \exp-\left(m(\lambda) + 2\int_0^\lambda g_q(s)\,ds\right), \qquad \lambda \geq 0.$$

By Lemma 3.6 and after some algebra, we get

$$c\lambda K'_{q,x}(\lambda) = -\exp-\left(\lambda x + m(\lambda) + \int_0^\lambda g_q(s)\,ds\right), \qquad \lambda > 0.$$

Next integrate the last equation to obtain

$$K_{q,x}(\lambda) = K_0 + \int_\lambda^\infty \frac{dt}{ct} e^{-tx - m(t) - \int_0^t g_q(s)\,ds}$$

$$= K_0 + \int_{\theta(\lambda)}^\infty dv\, \frac{\varphi'^2(v)}{c\varphi(v)} e^{-x\varphi(v) - \int_0^v w_q(u)\,du}, \qquad \lambda > 0,$$

where the last equality stems from the changes $u = \theta(s)$ and $v = \theta(t)$, and $K_0$ is an integration constant. Therefore, by the same kind of change, we get

$$F'_{q,x}(\lambda) = K_0 e^{m(\lambda) + 2\int_0^{\theta(\lambda)} w_q(z)\,dz}$$

$$+ e^{m(\lambda) + \int_0^{\theta(\lambda)} w_q(z)\,dz} \int_{\theta(\lambda)}^\infty ds\, \frac{\varphi'^2(s)}{c\varphi(s)} e^{-x\varphi(s) - \int_{\theta(\lambda)}^s w_q(z)\,dz}, \qquad \lambda > 0.$$

Check that this last function is integrable at 0. Indeed, $F'_{q,x}(\lambda) \sim K_{q,x}(\lambda)$ as $\lambda \searrow 0$, and the first expression we got for $K_{q,x}$ provides the equivalent $K_{q,x}(\lambda) \sim c^{-1} \ln(\lambda)$. Then because $F_{q,x}(0) = G_{q,x}(0) = 1/q$, the integration of $F'_{q,x}$ reads

$$F_{q,x}(\lambda) = \frac{1}{q} + K_0 \int_0^{\theta(\lambda)} dt\, e^{2\int_0^t w_q(z)\,dz}$$

$$+ \int_0^{\theta(\lambda)} dt\, e^{\int_0^t w_q(z)\,dz} \int_z^\infty ds\, \frac{\varphi'^2(s)}{c\varphi(s)} e^{-x\varphi(s) - \int_t^s w_q(z)\,dz}, \qquad \lambda \geq 0.$$

We thus get the following expression for $G_{q,x}$:

$$G_{q,x}(\lambda) = q^{-1}e^{-\int_0^{\theta(\lambda)} w_q(z)\,dz} + K_0 e^{-\int_0^{\theta(\lambda)} w_q(z)\,dz} \int_0^{\theta(\lambda)} dt\, e^{2\int_0^t w_q(z)\,dz}$$

$$+ e^{-\int_0^{\theta(\lambda)} w_q(z)\,dz} \int_0^{\theta(\lambda)} dt\, e^{\int_0^t w_q(z)\,dz} \int_t^\infty ds\, r^2(s)\, e^{-x\varphi(s) - \int_t^s w_q(z)\,dz},$$

$$\lambda \geq 0.$$

Recall that $qr^2 = w_q^2 - w_q'$, so that an integration by parts yields

$$\int_0^{\theta(\lambda)} dt\, e^{\int_0^t w_q(z)\,dz} \int_t^\infty ds\, qr^2(s)\, e^{-\int_t^s w_q(z)\,dz}$$

(13)
$$= \int_0^{\theta(\lambda)} dt\, w_q(t) e^{\int_0^t w_q(z)\,dz}$$

$$= e^{\int_0^{\theta(\lambda)} w_q(z)\,dz} - 1.$$

Since $qG_{q,0}$ is constant equal to 1, the constant integration $K_0$ must be 0. In addition, the previous display yields for any nonnegative $\lambda$,

$$qG_{q,x}(\lambda) = e^{-\int_0^{\theta(\lambda)} w_q(z)\,dz}$$

$$+ \int_0^{\theta(\lambda)} dt\, e^{-\int_t^{\theta(\lambda)} w_q(z)\,dz}$$

(14)
$$\times \int_t^\infty ds\, qr^2(s) e^{-x\varphi(s)} e^{-\int_t^s w_q(z)\,dz}$$

$$= 1 - \int_0^{\theta(\lambda)} dt\, e^{-\int_t^{\theta(\lambda)} w_q(z)\,dz}$$

$$\times \int_t^\infty ds\, qr^2(s)(1 - e^{-x\varphi(s)}) e^{-\int_t^s w_q(z)\,dz},$$

which ends the proof of (7). For (8), recall that $T_a$ is a.s. finite upon the condition given in Theorem 3.5. Observe that

$$G_{q,x}(\lambda) = \int_0^\infty dt\, e^{-qt} \mathbb{E}_x(\mathbf{1}_{Z_t=0} + e^{-\lambda Z_t}\mathbf{1}_{Z_t>0}), \qquad \lambda \geq 0.$$

By a double application of Lebesgue's convergence theorem, this quantity converges to

$$\int_0^\infty dt\, e^{-qt} \mathbb{P}_x(Z_t = 0),$$

so that, integrating by parts the last quantity,

$$\lim_{\lambda \to \infty} qG_{q,x}(\lambda) = \mathbb{E}_x(\exp(-qT_a)).$$

For the expectation of $T_a$ (9), we have to handle the Laplace exponent $H(t, \lambda)$ of the marginal distribution of the LB-process. Fubini's theorem and Beppo Levi's theorem, indeed, yield

$$\lim_{\lambda \to \infty} \int_0^\infty (1 - H(t, \lambda)) \, dt = \mathbb{E} \int_0^\infty \mathbf{1}_{Z_t \neq 0} \, dt = \mathbb{E}(T_a).$$

It is clear (see proof of Lemma 3.6) that

$$\frac{\partial H}{\partial t} = -\psi(\lambda) \frac{\partial H}{\partial \lambda} + c\lambda \frac{\partial^2 H}{\partial \lambda^2}, \qquad t, \lambda \geq 0,$$

which yields, by integrating $\frac{\partial H}{\partial \lambda} e^{-m(\lambda)}$ w.r.t. $\lambda$,

$$\frac{\partial H}{\partial \lambda} = -e^{m(\lambda)} \int_\lambda^\infty \frac{ds}{cs} \frac{\partial H}{\partial t}(t, s) e^{-m(s)}, \qquad t, \lambda \geq 0.$$

Integrate this last equation w.r.t. $\lambda$ and then w.r.t. $t$, and recall that $H(t, 0) = 1$, to get, by Fubini's theorem,

$$\int_0^\infty (1 - H(t, \lambda)) \, dt = \int_0^\lambda du \, e^{m(u)} \int_u^\infty \frac{ds}{cs} e^{-m(s)} (1 - e^{-sx}), \qquad \lambda \geq 0.$$

The result follows by letting $\lambda$ go to $\infty$ and applying Lebesgue's convergence theorem. An elementary variable change yields (9). □

PROOF OF COROLLARY 3.10. Equations (10) and (11) of the corollary follow from the same arguments that supported the theorem.

Let us prove the convergence of $\mathbb{E}_x(\exp(-\lambda Z_\tau))$ to (12). Observe that if $\tau$ is an independent exponential r.v. with parameter $q$, then $\mathbb{E}_x(\exp(-\lambda Z_\tau)) = qG_{q,x}(\lambda)$. The convergence then stems from (14) and the integration by parts used in (13). Namely,

$$\begin{aligned}
0 &\leq \mathbb{E}_x(\exp(-\lambda Z_\tau)) - e^{-\int_0^{\theta(\lambda)} w_q(z) \, dz} \\
&\leq \int_0^{\theta(\lambda)} dt \, e^{-\int_t^{\theta(\lambda)} w_q(z) \, dz} e^{-x\varphi(t)} \int_t^\infty ds \, qr^2(s) e^{-\int_t^s w_q(z) \, dz} \\
&= \int_0^{\theta(\lambda)} dt \, e^{-\int_t^{\theta(\lambda)} w_q(z) \, dz} e^{-x\varphi(t)} w_q(t).
\end{aligned}$$

Now the function

$$t \mapsto w_q(t) e^{-\int_t^{\theta(\lambda)} w_q(z) \, dz}$$

is integrable on $[0, \theta(\lambda)]$, so that, by Lebesgue's convergence theorem, the last displayed upper bound converges to 0 as $x \to \infty$.

Finally, the existence of (and weak convergence to) the LB-process starting from infinity stem from the previous convergence of Laplace transforms and standard theory [8]. □

## 4. Proofs.

4.1. *Proofs of Lemma* 3.3, *Theorems* 3.4 *and* 3.5.

PROOF OF LEMMA 3.3. First observe that

$$-m(\lambda) = \int_0^\lambda \frac{ds}{cs}\left(\delta s + \int_0^\infty \Pi(dr)(1-e^{-rs})\right)$$

$$= \frac{\delta\lambda}{c} + \int_0^\infty \Pi(dr)\int_0^{\lambda r} du\, \frac{1-e^{-u}}{cu}$$

$$= \frac{\delta}{c}\lambda + \int_0^\infty (1-e^{-\lambda r})\frac{\bar\Pi(r)}{cr}\,dr.$$

Note that the positive measure $\Lambda$ defined by $\Lambda(dr) = (cr)^{-1}\bar\Pi(r)\,dr$, $r > 0$, is, under condition (L), a Lévy measure. Indeed,

$$\int_1^\infty \Lambda(dr) = c^{-1}\int_1^\infty \ln(u)\Pi(du) < \infty,$$

$$\int_0^1 r\Lambda(dr) = c^{-1}\int_0^\infty (1\wedge u)\Pi(du) < \infty,$$

so that one can define $Y$ as the subordinator with Laplace exponent $-m$ ([5], Chapter III). As a consequence, the probability measure $\nu$ is well defined as the law of $Y_1$ (so it is infinitely divisible), and the measure $\Lambda$ always has infinite mass so that $Y_1 > \delta/c$ a.s. $\square$

PROOF OF THEOREM 3.4. To see that

(15) $$(\partial) \Leftrightarrow \int_0^\infty r^{-1}\nu(dr) < \infty,$$

we use the following equality (where both sides can be infinite):

$$\int_0^\infty r^{-1}\nu(dr) = \int_0^\infty e^{m(\lambda)}\,d\lambda,$$

so that (15) reduces to $(\partial) \Leftrightarrow \int^\infty d\lambda\, \exp(m(\lambda)) < \infty$.

Recall that $m(\lambda) = \int_0^\lambda (cs)^{-1}\psi(s)\,ds \leq 0$, $\lambda \geq 0$. It is clear from Lemma 3.3 that whenever $\delta \neq 0$, $\int^\infty d\lambda\, \exp(m(\lambda))$ is finite.

Next assume that $\Pi$ has finite mass $\rho$. Since the case when $\Pi$ has infinite mass can easily be derived from the case $\rho > c$ by a truncation argument, it only remains to prove that when $\delta = 0$, $\int^\infty d\lambda\, \exp(m(\lambda)) < \infty$ iff $\rho > c$. Then recall Lemma 3.3 and write, for any $\lambda > 1$,

$$-m(\lambda) = \int_0^\infty (1-e^{-\lambda r})\frac{\bar\Pi(r)}{cr}\,dr.$$

Next pick $x > 0$, divide $(0, \infty)$ into $(0, x/\lambda]$, $(x/\lambda, x]$ and $(x, \infty)$, and change variables $(u = \lambda r)$, to get

$$-m(\lambda) - \int_{x/\lambda}^{x} \frac{\bar{\Pi}(r)}{cr} dr$$
$$= \int_{0}^{x} (1 - e^{-u}) \frac{\bar{\Pi}(u/\lambda)}{cu} du$$
$$- \int_{x}^{\lambda x} e^{-u} \frac{\bar{\Pi}(u/\lambda)}{cu} du + \int_{x}^{\infty} (1 - e^{-\lambda r}) \frac{\bar{\Pi}(r)}{cr} dr.$$

Since $\bar{\Pi}$ is positive and decreases from $\bar{\Pi}(0+) = \rho$, an appeal to Lebesgue's convergence theorem yields

$$\lim_{\lambda \to \infty} -m(\lambda) - \int_{x/\lambda}^{x} \frac{\bar{\Pi}(r)}{cr} dr$$
(16)
$$= \rho \int_{0}^{x} (1 - e^{-u}) \frac{du}{cu} - \rho \int_{x}^{\infty} e^{-u} \frac{du}{cu} + \int_{x}^{\infty} \frac{\bar{\Pi}(r)}{cr} dr$$
$$\doteq k(x),$$

and $k(x)$ is finite. Now the same properties of $\bar{\Pi}$ ensure that

(17) $$\frac{\bar{\Pi}(x)}{c} \ln(\lambda) \leq \int_{x/\lambda}^{x} \frac{\bar{\Pi}(r)}{cr} dr \leq \frac{\rho}{c} \ln(\lambda), \qquad x > 0, \lambda > 1.$$

Using (16) and the left-hand side inequality in (17), we get

$$\limsup_{\lambda \to \infty} e^{m(\lambda)} \lambda^{\bar{\Pi}(x)/c} \leq \exp(-k(x)) < \infty.$$

Assume first that $\rho > c$ and pick $x$ small enough to have $\bar{\Pi}(x)/c > 1$, then the last equation shows that $\int^{\infty} d\lambda \exp(m(\lambda)) < \infty$.

Using (16) and the right-hand side inequality in (17), we get

$$\liminf_{\lambda \to \infty} e^{m(\lambda)} \lambda^{\rho/c} \geq \exp(-k(x)) > 0.$$

If $\rho \leq c$, the last equation then shows that $\int^{\infty} d\lambda \exp(m(\lambda)) = \infty$. This ends the proof of (15).

We next point out that, ($\partial$) holding or not, $Z$ is recurrent in $(\delta/c, \infty)$. Indeed, it is known [34] that $R$ oscillates in $(\delta/c, \infty)$ and converges in distribution to $\nu$, so that $T_0 = \infty$ a.s., $\lim_{t \to \infty} \eta_t = \infty$ a.s., and $C$ is defined on the entire $[0, \infty)$, so $Z = R \circ C$ also oscillates in $(\delta/c, \infty)$.

Let us characterize the invariant measures. According to Definition 3.2, the infinitesimal generator $Q$ of $Z$ satisfies, for any positive $\lambda$,

$$Qe_\lambda(z) = zAe_\lambda(z) - cz^2 e'_\lambda(z) = (\psi(\lambda) + c\lambda z)ze_\lambda(z), \qquad z \geq 0,$$

where $e_\lambda(z) = e^{-\lambda z}$. Now for any positive measure $\zeta$ on $[0,\infty)$, let $\chi$ denote the Laplace transform of $r\zeta(dr)$. Then $\zeta$ is invariant for $Z$ iff $\zeta Q = 0$ and

$$\begin{aligned}
\zeta Q = 0 &\iff \zeta Q e_\lambda = 0 \quad \forall \lambda > 0 \\
&\iff \int_0^\infty \zeta(dz)(\psi(\lambda) + c\lambda z) z e_\lambda(z) = 0 \quad \forall \lambda > 0 \\
&\iff \psi(\lambda)\chi(\lambda) - c\lambda\chi'(\lambda) = 0, \quad \lambda > 0.
\end{aligned}$$

The resolution of this ordinary differential equation proves that $\chi(\lambda) = k\exp(m(\lambda))$ for some constant $k$, so that

(18) $$r\zeta(dr) = k\nu(dr), \quad r > 0,$$

with a possible Dirac mass at 0 for $\zeta$. We can now prove the dichotomy between (i) (convergence in distribution) and (ii) (convergence to 0 in probability).

(i) Thanks to (15), since $(\partial)$ holds, we know that $\int_0^\infty r^{-1}\nu(dr) < \infty$. Equation (18) then ensures that $\mu$ is the unique invariant probability measure.

(ii) In the case when $(\partial)$ does not hold, $\int_0^\infty r^{-1}\nu(dr) = \infty$ and all invariant positive measures on $(0,\infty)$ are nonintegrable at $0+$, so that $\lim_{t\to\infty} \mathbb{P}(Z_t > \epsilon) = 0$ for any positive $\epsilon$.

Recall that $Z$ oscillates in $(0,\infty)$ to conclude it is null-recurrent. $\square$

PROOF OF THEOREM 3.5. In the case when $X$ is not a subordinator, the Ornstein–Uhlenbeck type process $R$ still oscillates [34] with stationary measure $\nu$, but this time $\nu$ charges the negative half-line so that $\liminf_{t\to\infty} R_t < 0$ a.s., and $T_0 < \infty$ a.s. It is thus clear that $\lim_{t\to\infty} Z_t = 0$ a.s., the absorption depending on whether the integral $\int_0^{T_0} ds/R_s (= T_a)$ converges or not. To check the criterion for absorption given in the theorem, we make use of the function $G_{q,x}$. First define the positive function $g_{q,x}$ as

$$g_{q,x}(\lambda) = -\frac{G'_{q,x}}{G_{q,x}}(\lambda), \quad \lambda \geq 0.$$

In the rest of the proof, we will drop the indices whenever it is not ambiguous to do so. Recall that

$$\mathbb{E}(\exp(-qT_a)) = \lim_{\lambda\to\infty} qG(\lambda).$$

Since $qG(\lambda) = \exp(-\int_0^\lambda g(s)\,ds)$, it suffices to show that

$$\int^\infty g(\lambda)\,d\lambda \quad \text{and} \quad \int^\infty d\lambda/\psi(\lambda) \quad \text{have the same nature.}$$

Thanks to Lemma 3.6, it is not difficult to prove that $g$ solves

$$c\lambda(-y' + y^2) + \psi y = q - \frac{e^{-\lambda x}}{G(\lambda)}, \qquad \lambda > 0.$$

We first show that $G$ decreases more slowly than any negative exponential function of parameter, say, $z$

$$\lim_{\lambda \to \infty} \frac{e^{-z\lambda}}{G(\lambda)} = 0.$$

To this aim, define

$$\sigma = \sup\{t > 0 : Z_t > z/2\},$$

which is finite a.s. since $Z$ goes to 0. Next observe that

$$G(\lambda) \geq \int_0^\infty dt\, e^{-qt}\, \mathbb{E}(e^{-\lambda Z_t}, t > \sigma)$$

$$\geq e^{-\lambda z/2} \int_0^\infty dt\, e^{-qt} \mathbb{P}(t > \sigma)$$

$$= e^{-\lambda z/2} q^{-1} \mathbb{E}(e^{-q\sigma}),$$

which proves our claim. Since $X$ is not a subordinator, $\lambda/\psi(\lambda)$ is bounded from above (and positive for sufficiently large $\lambda$), so that

$$h(\lambda) \doteq \frac{e^{-\lambda x}}{\psi(\lambda) G(\lambda)} = o(e^{-\lambda x/2}) \qquad \text{as } \lambda \to \infty.$$

We will use this last comparison along with the following rearrangement of the differential equation that $g$ solves

(19) $\qquad \psi(\lambda)(g + h)(\lambda) - q - c\lambda g'(\lambda) = -c\lambda g^2(\lambda), \qquad \lambda > 0.$

On the vector space $\mathcal{C}_P$ of continuous real functions defined on $[0, \infty)$ with (at most) polynomial growth at $\infty$, define the scalar product $\langle \cdot, \cdot \rangle$ by

$$\langle u, v \rangle = \frac{\int_0^\infty dt\, e^{-qt} \mathbb{E}_x(u(Z_t) v(Z_t) e^{-\lambda Z_t})}{G_{q,x}(\lambda)}, \qquad u, v \in \mathcal{C}_P.$$

In particular, if $\bar{1}$ denotes the constant function equal to 1 and $I$ the identity function on $[0, \infty)$, then $\langle \bar{1}, \bar{1} \rangle = 1$ and

$$g'(\lambda) = -\frac{G''}{G}(\lambda) + \left(\frac{G'}{G}(\lambda)\right)^2 = -\langle \bar{1}, \bar{1} \rangle \langle I, I \rangle + \langle I, \bar{1} \rangle^2 \leq 0,$$

by the Cauchy–Schwarz inequality. Thus, the positive function $g$ is decreasing on $[0, \infty)$. If $l \geq 0$ is its limit at $\infty$, since $\limsup_{\lambda \to \infty} g'(\lambda) = 0$ and $\lim_{\lambda \to \infty} \lambda/\psi(\lambda) \in [0, \infty)$, (19) yields $l \leq 0$, so that $l = 0$.

Thanks to this same equation, and because $l = 0$, note that

$$\lim_{\lambda \to \infty} g(\lambda)^{-1}(g(\lambda) + h(\lambda) - \psi(\lambda)^{-1}(q + c\lambda g'(\lambda))) = -\lim_{\lambda \to \infty} \frac{c\lambda}{\psi(\lambda)} g(\lambda) = 0.$$

As a consequence,

$$\int^\infty g(\lambda) \, d\lambda \quad \text{and}$$

$$\int^\infty (\psi(\lambda)^{-1}(q + c\lambda g'(\lambda)) - h(\lambda)) \, d\lambda \qquad \text{have the same nature.}$$

But $h(\lambda) = o(e^{-\lambda x/2})$, $\lim_{\lambda \to \infty} \lambda/\psi(\lambda) \in [0, \infty)$ and $\lim_{\lambda \to \infty} \downarrow g(\lambda) = l = 0$, so that

$$\int^\infty \lambda \psi(\lambda)^{-1} g'(\lambda) \, d\lambda \quad \text{and} \quad \int^\infty h(\lambda) \, d\lambda \qquad \text{both converge,}$$

which reduces the comparison to $\int^\infty d\lambda/\psi(\lambda)$. $\square$

4.2. *Proofs of Lemmas* 3.6, 3.7 *and* 2.1.

PROOF OF LEMMA 3.6. Recall (6), and apply it to the function $e_\lambda$ defined as $e_\lambda(z) = e^{-\lambda z}$. It is well known that $e_\lambda$ is in the domain of $A$ and that, for any positive $\lambda$,

$$Ae_\lambda(z) = \psi(\lambda) e_\lambda(z), \qquad z \geq 0.$$

As a consequence,

$$qG_{q,x}(\lambda) = e_\lambda(x) + \int_0^\infty dt \, e^{-qt} \mathbb{E}_x(\psi(\lambda) Z_t e^{-\lambda Z_t} + c\lambda Z_t^2 e^{-\lambda Z_t})$$
$$= e_\lambda(x) - \psi(\lambda) G'_{q,x}(\lambda) + c\lambda G''_{q,x}(\lambda),$$

which concludes the proof. $\square$

PROOF OF LEMMA 3.7. Since $f_q$ solves the homogeneous equation $(E_q^{(h)})$ associated to $(E_q)$, we have

$$-c\lambda f''_q(\lambda) + \psi(\lambda) f'_q(\lambda) + q f_q(\lambda) = 0, \qquad \lambda \in I.$$

It is an old trick to consider $g_q = -f'_q/f_q$, for then $f''_q/f = -g'_q + (f'_q/f_q)^2$. For $\lambda \in I_0$ [$f_q(\lambda) \neq 0$], divide the previously displayed equality by $-c\lambda f_q(\lambda)$, and get

$$-g'_q(\lambda) + g_q^2(\lambda) + \frac{\psi(\lambda)}{c\lambda} g_q(\lambda) = \frac{q}{c\lambda}, \qquad \lambda \in I_0.$$

For any $g_q$ satisfying the last equation on some open interval $J$, we consider the change of function $h_q(\lambda) = e^{-m \circ \varphi(\lambda)} g_q \circ \varphi(\lambda)$, $\lambda \in \theta(J)$, so that

$$h_q'(\lambda) = e^{-2m \circ \varphi(\lambda)}(-m' \circ \varphi(\lambda) g_q \circ \varphi(\lambda) + g_q' \circ \varphi(\lambda))$$

$$= e^{-2m \circ \varphi(\lambda)}\left(g_q^2 \circ \varphi(\lambda) - \frac{q}{c\varphi(\lambda)}\right)$$

$$= h_q^2(\lambda) - \frac{q\varphi'(\lambda)^2}{c\varphi(\lambda)},$$

which indeed proves that $h_q$ solves $(\epsilon_q')$ on $\theta(J)$. □

PROOF OF LEMMA 2.1.

NOTATION. For any real number $t$ and real functions $f, g$ which do not vanish in a pointed neighborhood of $t$, we will write $f(s) \sim g(s)$ [resp. $f(s) \asymp g(s)$] as $s \to t$ to mean

$$\lim_{s \to t} \frac{f(s)}{g(s)} = 1 \qquad (\text{resp.} = \text{some nonzero, finite real number}).$$

Recall that $(\epsilon_q')$ is the following Riccati differential equation:

$$(\epsilon_q') \qquad\qquad y' - y^2 = -qr^2,$$

where

$$r(s) = \frac{|\varphi'(s)|}{\sqrt{c\varphi(s)(1 - \varphi(s))}}, \qquad s \in (0, \xi),$$

and $\xi < \infty \Leftrightarrow d < c$. By the Cauchy–Lipschitz theorem (see, e.g., Corollaries 1 and 2, Chapter 7.2, page 93 in [1]), for any $T \in (0, \xi)$ and $a \in \mathbb{R}$, there is a unique maximal solution $(I, f)$ to $(\epsilon_q')$ with boundary condition $f(T) = a$. We denote by $f_{(T,a)}$ this maximal solution and by $I_{(T,a)}$ the corresponding maximal interval. We will use repeatedly, and sometimes implicitly, the fact that the orbits of different solutions never intersect, therefore inducing a total ordering on any set of solutions defined on the same interval.

Let us give the outline of the proof. We state three lemmas and then provide the proof of Lemma 2.1. The first lemma studies the graph of $r$ and shows, in particular, that $r$ is monotone on neighborhoods of $0+$ and $\xi-$, and that $\lim_{t \to 0+} r(t) = +\infty$. The second one is a technical result ensuring, in particular, that, for any $T \in (0, \xi)$, the maximal interval $I_{(T,0)}$ is the entire $(0, \xi)$. The third lemma asserts that the family $(f_{(T,0)})_T$ is uniformly increasing and bounded from above, allowing us to define $w_q$ as their increasing limit as $T \to \xi-$.

The rest of the section will then be devoted to proving that $w_q$ solves $(\epsilon'_q)$, and that it is the unique nonnegative solution vanishing at $\xi-$, and to studying its behavior at $0+$ and $\xi-$.

LEMMA 4.1. *The function $r$ is monotone on a neighborhood of $\xi-$ and monotone increasing on a neighborhood of $0+$. Furthermore,*

$$\lim_{s \to 0+} r(s) = +\infty,$$

$$\lim_{s \to \xi-} r(s) = \begin{cases} +\infty, & \text{if } d < c/2, \\ 1/\beta\sqrt{c}, & \text{if } d = c/2, \\ 0, & \text{if } d > c/2. \end{cases}$$

We point out that when $d \geq c$, one has $\xi = +\infty$, $\lim_{s \to +\infty} r(s) = 0$, and $\lim_{s \to 0+} r(s) = +\infty$, which is the same situation as in the continuous setting. The sequel of the proof in the case when $d \geq c$ is thus exactly the same in the continuous setting.

PROOF OF LEMMA 4.1. Set

$$H(s) \doteq r^2 \circ \theta(s) = \frac{e^{-2m(s)}}{cs(1-s)}, \qquad s \in (0,1),$$

and recall (3) to see that

$$H(s) \sim \frac{1}{c\beta^2} s^{2d/c-1} \quad \text{as } s \to 0+ \quad \text{and} \quad H(s) \sim \frac{1}{c(1-s)} \quad \text{as } s \to 1-.$$

Then differentiate $H$ to check that it is monotone on neighborhoods of $0$ and $\xi$, so that the same holds for $r$, by monotonicity of $\theta$. □

LEMMA 4.2. (i) *For any $T \in (0,\xi)$ and $a \in \mathbb{R}$,*

$$a \leq 0 \implies [T,\xi) \subset I_{(T,a)} \text{ and } f_{(T,a)} \text{ is negative on } [T,\xi),$$

$$a \geq 0 \implies (0,T] \subset I_{(T,a)} \text{ and } f_{(T,a)} \text{ is positive on } (0,T].$$

(ii) *If $d < c/2$, there is a neighborhood $\mathcal{V}$ of $\xi$ such that for any $T \in \mathcal{V}$ and $a \in \mathbb{R}$,*

$$a \leq \sqrt{q} r(T) \implies [T,\xi) \subset I_{(T,a)} \quad \text{and} \quad \lim_{t \to \xi-} f_{(T,a)}(t) \text{ exists and is finite.}$$

PROOF. First note that for any $T \in (0,\xi)$, $f_{(T,0)}(T) = 0$ and $f'_{(T,0)}(T) = -qr^2(T) < 0$. Consequently, a maximal solution of $(\epsilon'_q)$ vanishes, at most, once on its interval of definition.

(i) Let $a \leq 0$. We just saw why $f \doteq f_{(T,a)}$ remains negative for $t > T$. Now assume that the supremum $\sigma$ of $I_{(T,a)}$ is in $(T, \xi)$. Since $f$ is a maximal solution and is negative on $(T, \xi)$, we deduce that $\lim_{t \to \sigma-} f(t) = -\infty$. But this would mean that $f$ is below $-\sqrt{q}r$ on some neighborhood of $\sigma$, and thus increasing. This shows that $\sigma = \xi$. The same arguments hold for the other implication $(a \geq 0)$.

(ii) Assume $a \in [-\sqrt{q}r(T), \sqrt{q}r(T)]$. When $d < c/2$ (Lemma 4.1), $r$ increases on some neighborhood $\mathcal{V}$ of $\xi < \infty$, and goes to $\infty$ at $\xi-$. For this reason, $f \doteq f_{(T,a)}$ decreases $(f' = f^2 - qr^2)$, so it remains below $\sqrt{q}r$, but never reaches $-\sqrt{q}r$, since at a hitting time of $-\sqrt{q}r$, $f'$ would vanish, whereas $-\sqrt{q}r$ decreases. So $f$ is defined, as claimed, at least on $[T, \xi)$ even if $a > 0$, and for any $t$ in a neighborhood of $\xi$, $f^2(t) < qr^2(t)$. The same holds when $a < -\sqrt{q}r(T)$, for then $f$ increases until it reaches $-\sqrt{q}r$.

To get the existence of a finite limit at $\xi-$, it is thus equivalent, by $f' - f^2 = -qr^2$ and $|f'| < qr^2$, to prove that $\int^{\xi} r^2 < \infty$. To see this, use (3) and get, as $s \to \xi-$,

$$\varphi(s) \asymp (\xi - s)^{1/(1-d/c)} \qquad \text{so that } r^2(s) \asymp (\xi - s)^{-(1-2d/c)/(1-d/c)}.$$

The result follows by observing that if $d \in (0, c/2)$, then $-\frac{1-2d/c}{1-d/c} \in (-1, 0)$.
□

As a consequence of the previous lemma, for any $T \in (0, \xi)$, $f_{(T,0)}$ is defined on the entire $(0, \xi)$, that is, $I_{(T,0)} = (0, \xi)$. Next set

$$\kappa_T \doteq f_{(T,0)}, \qquad T \in (0, \xi).$$

LEMMA 4.3. *The sequence* $(\kappa_T)_T$ *is uniformly increasing and uniformly bounded from above on* $(0, \xi)$. *One can then define the function* $w_q$ *as the uniformly increasing limit*

$$w_q = \lim_{T \uparrow \xi} \uparrow \kappa_T.$$

*The function* $w_q$ *is positive on* $(0, \xi)$ *and* $\lim_{t \to \xi-} w_q(t) = 0$.

PROOF. The orbits $(\kappa_T)_T$ are totally ordered because they never intersect, and increase as $T$ increases because $\kappa'_T(T) = -qr^2(T) < 0$.

We next intend to prove that the family $(\kappa_T)_T$ is uniformly bounded from above.

(i) $d \geq c$. In this case, $\xi = \infty$, and there is a $t_0$ such that $r$ decreases on $[t_0, \infty)$. Thanks to Lemma 4.2(i), $h_0 \doteq f_{(t_0, \sqrt{q}r(t_0))}$ is defined (at least) on $(0, t_0]$ and is positive. Since different orbits cannot intersect, for any $T \in (0, \infty)$, $\kappa_T$ is below $h_0$ on $(0, t_0]$. Conclude by noticing that on $[t_0, \infty)$, $\kappa_T < \sqrt{q}r$ (otherwise $\kappa_T$ would increase on $[t_0, \infty)$, which contradicts Lemma 4.2(i)).

(ii) $d \in [c/2, c)$. Recall Lemma 4.1 to see that $r$ can be extended continuously to $(0, \xi]$. In the same vein as in Lemma 4.2(i), the maximal solution $\kappa_\xi$ of $(\epsilon'_q)$ vanishing at $\xi$ is defined on $(0, \xi]$. It is then obvious that $\kappa_T$ is below $\kappa_\xi$ for any $T < \xi$.

(iii) $d < c/2$. Thanks to Lemma 4.2(ii), there is a $t_0$ close enough from $\xi$ such that $h_0 \doteq f_{(t_0, \sqrt{q}r(t_0))}$ has a positive limit at $\xi-$. This function $h_0$ is defined on $(0, \xi)$ and since it never vanishes, it never intersects any $\kappa_T$. As a consequence, the family $(\kappa_T)_T$ is bounded from above by $h_0$.

Now for any $t \in (0, \xi)$ and $T > t$, $\kappa_T(t) > 0$ [Lemma 4.2(i)]. Since $w_q(t)$ is the increasing limit of $\kappa_T(t)$ as $T$ increases, $w_q(t) > 0$. The function $w_q$ vanishes at $\xi-$ by continuity of the flow. $\square$

We hereafter conclude the proof of Lemma 2.1. We first prove that for any positive $q$, $w_q$ is the unique nonnegative solution of $(\epsilon'_q)$ defined on $(0, \xi)$ and vanishing at $\xi-$.

Let us show the uniqueness of such a solution. Assume $h$ and $k$ are two nonnegative solutions of $(\epsilon'_q)$ on a neighborhood of $\xi$, who both vanish at $\xi-$. If there is some $t_0$ such that $h(t_0) \neq k(t_0)$, say $h(t_0) > k(t_0)$, then $u \doteq h - k$ is positive on $(t_0, \xi)$ and so is $v \doteq h + k$. Since $h$ and $k$ are both solutions of $(\epsilon'_q)$, $u$ and $v$ are such that $u' = uv$, so that $u$ is positive and increasing on $(t_0, \xi)$. This contradicts the fact that $\lim_{t \to \xi-} u(t) = 0$.

To show that $w_q$ is the solution of $(\epsilon'_q)$, fix $\varepsilon$ in $(0, \xi)$, and set $\eta \doteq w_q(\varepsilon)$ and $h_\varepsilon \doteq f_{(\varepsilon, \eta)}$. Next suppose that for some $t_0 \in I_{(\varepsilon, \eta)}$, $t_0 > \varepsilon$ and $h_\varepsilon(t_0) \neq w_q(t_0)$. Since $w_q$ is by construction the supremum of a family of orbits indexed by a continuous set, it is necessary that $h_\varepsilon(t_0) > w_q(t_0)$. Pick then some $a \in (w_q(t_0), h_\varepsilon(t_0))$. We know from Lemma 4.2(i) that $f \doteq f_{(t_0, a)}$ is defined at least on $(0, t_0]$. Since $f(t_0) > w_q(t_0)$, $f(t_0) > \kappa_T(t_0)$ for all $T \in (0, \xi)$. This last inequality certainly holds for all $t \in (0, t_0]$, and thus, in particular, $f(\varepsilon) \geq w_q(\varepsilon) = \eta$. Reasoning in the same fashion with $h_\varepsilon$, we see that $f(\varepsilon) < h_\varepsilon(\varepsilon)$. But $h_\varepsilon(\varepsilon) = \eta$, which brings the contradiction ($\eta \leq f(\varepsilon) < \eta$), and proves that $w_q(t) = h_\varepsilon(t)$ for all $t \in I_{(\varepsilon, \eta)}$, $t \geq \varepsilon$. Now by the same reasoning as in the proof of Lemma 4.2, $h_\varepsilon$ is defined (at least) on $[\varepsilon, \xi)$, since otherwise it would go to $+\infty$ in finite time and there would, indeed, be some $t_0$ such that $h_\varepsilon(t_0) > w_q(t_0)$. As a consequence, $w_q = h_\varepsilon$ on $[\varepsilon, \xi)$. This ends the proof, since $\varepsilon$ is arbitrarily small.

We complete the proof of Lemma 2.1 by studying the behavior of $w_q$ at $0+$ and $\xi-$.

First note that because $w'_q = w_q^2 - qr^2$, $w_q$ increases when it is above $\sqrt{q}r$, and that $w'_q$ vanishes when both graphs intersect. Recall $r$ is positive and monotone on neighborhoods of 0 and $\xi$. As a consequence, the graph of $w_q$ may intersect that of $\sqrt{q}r$ at most once on each of these two neighborhoods. On the other hand, it cannot remain above $\sqrt{q}r$ on a whole neighborhood:

(i) of $0+$, because it would be increasing but tend to $\infty$ as $t \searrow 0+$;
(ii) of $\xi-$, because it would be increasing and tend to 0 as $t \nearrow \xi-$.

Now since $w_q < \sqrt{q}r$ on neighborhoods of $0+$ and $\xi-$, it decreases on these neighborhoods, and $\int_0^\xi w_q$ converges because $\int_0^\xi r$ does so. Indeed,

$$r(s) \sim \frac{-\varphi'(s)}{\sqrt{c(1-\varphi(s))}} \quad \text{as } s \to 0+ \quad \text{and} \quad r(s) \sim \frac{-\varphi'(s)}{\sqrt{c\varphi(s)}} \quad \text{as } s \to \xi-,$$

which are both integrable on their respective neighborhoods, since $\varphi$ has finite limits at $0+$ and $\xi-$. $\square$

## REFERENCES

[1] ARNOLD, V. I. (1992). *Ordinary Differential Equations*. Springer, Berlin. MR1162307
[2] ATHREYA, K. B. and NEY, P. E. (1972). *Branching Processes*. Springer, New York. MR373040
[3] ATHREYA, S. R. and SWART, J. M. (2005). Branching–coalescing particle systems. *Probab. Theory Related Fields*. To appear. Available at http://springerlink.metapress.com/ link.asp?id=hm8np9hw554h1jbt. MR2123250
[4] BERESTYCKI, J. (2004). Exchangeable fragmentation–coalescence processes and their equilibrium measures. *Electron. J. Probab.* **9** 770–824. Available at http://hal.ccsd.cnrs.fr/4ccsd-00001256. MR2110018
[5] BERTOIN, J. (1996). *Lévy Processes*. Cambridge Univ. Press.
[6] BERTOIN, J. and LE GALL, J. F. (2000). The Bolthausen–Sznitman coalescent and the genealogy of continuous-state branching processes. *Probab. Theory Related Fields* **117** 249–266. MR1771663
[7] CAMPBELL, R. B. (2003). A logistic branching process for population genetics. *J. Theor. Biol.* **225** 195–203. MR2077387
[8] CHUNG, K. L. (1970). *Lectures on Boundary Theory for Markov Chains.* Princeton Univ. Press. MR267644
[9] DIACONIS, P., MAYER-WOLF, E., ZEITOUNI, O. and ZERNER, M. P. W. (2004). The Poisson–Dirichlet law is the unique invariant distribution for uniform split-merge transformations. *Ann. Probab.* **32** 915–938. MR2044670
[10] DUQUESNE, T. and LE GALL, J.-F. (2002). Random trees, Lévy processes and spatial branching processes. *Astérisque* **281**. MR1954248
[11] DUSHOFF, J. (2000). Carrying capacity and demographic stochasticity: scaling behavior of the stochastic logistic model. *Theoret. Population Biology* **57** 59–65.
[12] DYNKIN, E. B. (1965). *Markov Processes.* Springer, Berlin. MR193671
[13] ETHERIDGE, A. M. (2004). Survival and extinction in a locally regulated population. *Ann. Appl. Probab.* **14** 188–214. MR2023020
[14] GOSSELIN, F. (2001). Asymptotic behavior of absorbing Markov chains conditional on nonabsorption for applications in conservation biology. *Ann. Appl. Probab.* **11** 261–284. MR1825466

[15] Grey, D. R. (1974). Asymptotic behaviour of continuous-time, continuous statespace branching processes. *J. Appl. Probab.* **11** 669–677. MR408016

[16] Hgns, G. (1997). On the quasi-stationary distribution of a stochastic Ricker model. *Stochastic Process. Appl.* **70** 243–263.

[17] Jagers, P. (1992). Stabilities and instabilities in population dynamics. *J. Appl. Probab.* **29** 770–780. MR1188534

[18] Jagers, P. (1997). Towards dependence in general branching processes. In *Classical and Modern Branching Processes*. Springer, New York. MR1601717

[19] Jirina, M. (1958). Stochastic branching processes with continuous state space. *Czechoslovak Math. J.* **8** 292–312. MR101554

[20] Joffe, A. and Métivier, M. (1986). Weak convergence of sequences of semimartingales with applications to multitype branching processes. *Adv. in Appl. Probab.* **18** 20–65.

[21] Keller, G., Kersting, G. and Rösler, U. (1987). On the asymptotic behaviour of discrete time stochastic growth processes. *Ann. Probab.* **15** 305–343. MR877606

[22] Kingman, J. F. C. (1982). The coalescent. *Stochastic Process. Appl.* **13** 235–248. MR671034

[23] Klebaner, F. C. (1997). Population and density dependent branching processes. In *Classical and Modern Branching Processes*. Springer, New York. MR1601729

[24] Lambert, A. (2001). The branching process conditioned to be never extinct. Ph.D. thesis, Univ. P. et M. Curie. Available at http://www.biologie.ens.fr/ecologie/coevolution/lambert/lambert/qprocess.pdf.

[25] Lambert, A. (2002). The genealogy of continuous-state branching processes with immigration. *Probab. Theory Related Fields* **122** 42–70. MR1883717

[26] Lambert, A. (2003). Coalescence times for the branching process. *Adv. in Appl. Probab.* **35** 1071–1089. MR2014270

[27] Lamperti, J. (1967). Continuous-state branching processes. *Bull. Amer. Math. Soc.* **73** 382–386. MR208685

[28] Le Gall, J. F. and Le Jan, Y. (1998). Branching processes in Lévy processes: The exploration process. *Ann. Probab.* **26** 213–252.

[29] Matis, J. H. and Kiffe, T. R. (2004). On stochastic logistic population growth models with immigration and multiple births. *Theoret. Population Biology* **65** 89–104.

[30] Nåsell, I. (2001). Extinction and quasi-stationarity in the Verhulst logistic model. *J. Theor. Biol.* **211** 11–27.

[31] O'Connell, N. (1995). The genealogy of branching processes and the age of our most recent common ancestor. *Adv. in Appl. Probab.* **27** 418–442. MR1334822

[32] Pitman, J. (2002). Poisson–Dirichlet and GEM invariant distributions for split-and-merge transformation of an interval partition. *Combin. Probab. Comput.* **11** 501–514. MR1930355

[33] Ramanan, K. and Zeitouni, O. (1999). The quasi-stationary distribution for small random perturbations of certain one-dimensional maps. *Stochastic Process. Appl.* **84** 25–51. MR1720096

[34] Sato, K. I. (1999). *Lévy Processes and Infinitely Divisible Distributions*. Cambridge Univ. Press.

[35] Wofsy, C. L. (1980). Behavior of limiting diffusions for density-dependent branching processes. *Biological Growth and Spread. Lecture Notes in Biomath.* **38** 130–137. Springer, Berlin. MR609353

Unit of Mathematical Evolutionary Biology  
Fonctionnement et Évolution  
  des Systmes Écologiques UMR 7625  
École Normale Supérieure  
46 rue d'Ulm  
F-75230 Paris Cedex 05  
France  
e-mail: amaury.lambert@ens.fr  
url: http://www.biologie.ens.fr/ecologie/ecoevolution/lambert/index.en.html